\documentstyle[verbatim,amstex,12pt,leqno]{article}

\begin{document}
\bibliographystyle{amsplain}

\begin{center}
	{\LARGE{}Using sums of squares to prove that certain 
	entire functions have only real zeros}\\[.3cm]
	by\\[.3cm]
	{\large{}George  Gasper\footnotemark}

\footnotetext{Supported in part by the National Science Foundation under
grant DMS-9103177. \par
 }
        
\vskip.3cm
{\it{Dedicated to the memory of Ralph P. Boas, Jr. (1912--1992)}}
\vskip.3cm
        {(May 17, 1993 version)}\\[.4cm]
\vskip.3cm    
\end{center}

\begin{abstract}
It is shown how sums of squares of real valued functions can be used to
give new proofs of the reality of the zeros of the Bessel functions
$J_\alpha (z)$ when $\alpha \ge -1,$ confluent hypergeometric functions
${}_0F_1(c\/; z)$ when $c>0$ or $0>c>-1$, Laguerre polynomials
$L_n^\alpha(z)$ when $\alpha \ge -2,$ and Jacobi polynomials
$P_n^{(\alpha,\beta)}(z)$ when $\alpha \ge -1$ and $ \beta \ge -1.$ Besides
yielding new inequalities for $|F(z)|^2,$ where $F(z)$ is one of these
functions, the derived identities lead to inequalities for $\partial
|F(z)|^2/\partial y$ and $\partial ^2 |F(z)|^2/\partial y^2,$ which
also give new proofs of the reality of the zeros.
\end{abstract}

\section{Introduction}

In a 1975 survey paper \cite{gg75} on positivity and special functions
it was shown how sums of squares of special functions could be used to
prove the nonnegativity of the Fej\'er kernel, the positivity of
integrals of Bessel functions \cite{gg75B} and of the Cotes' numbers
for some Jacobi abscissas, a Tur\'an type inequality for Bessel
functions, the Askey-Gasper inequality (cf. \cite{ag76}, \cite{ag86},
\cite{gg86}, \cite{gg89})
\begin{equation}\label{1.1}
	\sum_{k=0}^n P_k^{(\alpha,0)}(x)
	\ge 0,\quad \alpha> -2,\quad  -1\le x\le 1,
\end{equation}
which  de Branges \cite{de} employed to complete his proof of the Bieberbach
conjecture, and to prove the more general inequalities \cite{gg77}
\begin{equation}\label{1.2}
	\sum_{k=0}^n {(\lambda +1)_k\over k!}{(\lambda +1)_{n-k}\over (n-k)!}
	{P_k^{(\alpha,\beta)}(x)\over P_k^{(\beta,\alpha)}(1)}
	\ge 0, \quad -1\le x\le 1,
\end{equation} 
when $0\le \lambda\le \alpha +\beta$ and $ \beta \ge -1/2$.  
It was also pointed out in \cite{gg75} that, since one of
Jensen's necessary and sufficient conditions for the Riemann Hypothesis  to
hold  (given in P\'olya \cite{pol}) is the condition that 
\begin{equation}\label{1.3}
	\int^\infty_{-\infty}\int^\infty_{-\infty} \Phi(s)\Phi(t)
	e^{i(s+t)x} (s- t)^{2n}dsdt\ge 0, \quad -\infty<x<\infty,
\end{equation}
 \noindent
for $n = 0, 1, 2, 3, \dots ,$ where
\begin{equation}\label{1.4}
	\Phi(t)=2\sum^\infty_{k=1} (2k^4\pi^2 e^{9t/2} - 3k^2\pi
	e^{5t/2}) e^{-k^2\pi e^{2t}},
\end{equation}
 \noindent
and the above integral is a square when  $n = 0$, the method  of sums of
squares is suggested for proving (1.3).  

Another of Jensen's necessary and 
sufficient conditions for the Riemann Hypothesis to hold is that
\begin{equation}\label{1.5}
	\int^\infty_{-\infty}\int^\infty_{-\infty} 
	\Phi(s)\Phi(t) e^{i(s+t)x} e^{(s- t)y} (s-t)^2ds dt\ge0,
	\quad -\infty<x,y<\infty,
\end{equation}
which can also be written in the equivalent
form 
\begin{equation}\label{1.6}
	{\partial^2 \over \partial y^2}|\Xi(x+iy)|^2 \ge 0, 
	\quad -\infty<x,y<\infty, 
\end{equation}
with
\begin{equation}\label{1.7}
	\Xi(z)=\int_{-\infty}^{\infty} \Phi(t) \exp(izt) dt =
	2\int_0^{\infty} \Phi(t) \cos(zt) dt.
\end{equation}
That (1.6) is a sufficient condition for the Riemann $\Xi(z)$ function
to have only real zeros follows directly from observation that, since
$|\Xi(x+iy)|^2= \Xi(x+iy)\Xi(x-iy)$ is a nonnegative even function of
y, (1.6) implies that $|\Xi(x+iy)|^2$ is a nonnegative even convex
function of $y$ with its unique minimum at $y=0,$ and hence $\Xi (x+iy)
\ne 0$ whenever $y \ne 0$. If the function $\Phi(t)$ in (1.3) and (1.5)
is replaced by a function $\Psi(t)$ such that the conditions stated in
\cite[\S1]{pol}, are satisfied, then, by \cite[pp.~17,~18]{pol}, the
inequalities in (1.3) and (1.5) are necessary and sufficient conditions
for the Fourier (or cosine) transform of $\Psi(t)$ to have only real
zeros.  In 1913 Jensen \cite{jen} proved that each of the inequalities
\begin{equation}\label{1.8}
	\quad y{\partial \over \partial y }|F(x+i y)|^2 \ge 0, \quad
	{\partial^2 \over \partial y^2}|F(x+i y)|^2
	\ge 0,\quad -\infty<x, y<\infty,
\end{equation}
is necessary and sufficient for a real entire function $F(z) \not\equiv 0$ of 
genus 0 or 1 (cf. Boas \cite[Chapter~2]{bo}) to have only real zeros.
Also see Titchmarsh \cite{titch} and Varga \cite[Chapter~3]{var}.  

In view of these observations and the successes of the sums of squares
method (also see \cite{gg89a}, \cite[Chapter~8]{gr}),
since the early 1970's I have been investigating
how squares of real valued functions can 
be used to prove that certain entire functions have only real zeros and to 
prove inequalities of the form in (1.8).
In this paper I demonstrate how certain series expansions in sums of squares
of special functions give
new proofs of the reality of the zeros of the Bessel functions $J_\alpha (z)$
when $\alpha \ge -1,$ confluent hypergeometric functions ${}_0F_1(c\/; z)$
when $c>0$ or $0>c>-1$, Laguerre polynomials $L_n^\alpha(z)$ when 
$\alpha \ge -2,$ and Jacobi polynomials $P_n^{(\alpha,\beta)}(z)$
when $\alpha \ge -1$ and $ \beta \ge -1.$  Here,
as elsewhere, $z = x + i y$ is a complex variable and
$x$ and $y$ are real variables. For
the definitions of these functions and their properties, see
Erd\'elyi \cite{erd} and  Szeg\H o \cite{sz}. In addition, it will 
be shown that besides yielding new inequalities for $|F(z)|^2,$
 where $F(z)$ is one of these functions, the derived identities
lead to inequalities for $\partial |F(z)|^2/\partial y$ and 
$\partial ^2 |F(z)|^2/\partial y^2,$ which also give new proofs of the reality 
of the zeros.

\section{ Initial observations}
\setcounter{equation}{0}

In order to see how easily sums of squares
can be used to prove that all of the zeros of $\sin z$ and $\cos z$ are
real, it suffices to observe that we have the (easily verified) identities 
\begin{equation}\label{2.1}
	|\sin z|^2 = \sin^2 x + \sinh^2 y,
\end{equation}
\begin{equation}\label{2.2}
	|\cos z|^2 = \cos^2 x +\sinh^2 y
\end{equation} 
and to note that $\sinh y = (e^y - e^{-y})/2 > 0$ when
$y>0,$ and $\sinh y<0 $ when $y<0.$ 

One can also take  partial derivatives  of the identities
in (2.1) and (2.2) with respect to $y$ to obtain
\begin{equation}\label{2.3}
	{\partial \over \partial y }|\sin z|^2= 
{\partial \over \partial y}|\cos z|^2 = \sinh 2y 
\end{equation}
which shows that $|\sin z|^2$ and $|\cos z|^2$ are increasing (decreasing)
functions of $y$ when $y>0$ ($y<0),$ and to obtain 
\begin{equation}\label{2.4}
	{\partial^2 \over \partial y^2}|\sin z|^2 =
	{\partial^2 \over \partial y^2}|\cos z|^2 = 2\cosh 2y =
	2(\cosh^2 y +\sinh^2 y) \ge 2,
\end{equation}
which shows that  $|\sin z|^2$ and $|\cos z|^2$ are convex functions
of $y.$  Then, because  $|\sin z|^2$ and $|\cos z|^2$ are nonnegative
even functions of $y$, it immediately follows from (2.3) and (2.4) that
 $\sin z$ and $\cos z$ have only real zeros. 

Observe that the reality of the 
zeros of $\sin z$ and $\cos z$ also follows from the inequalities
\begin{equation}\label{2.5}
	|\sin z|^2 > \sin^2 x,\quad |\cos z|^2 > \cos^2 x,\quad (y \ne 0) 
\end{equation}
\begin{equation}\label{2.6}
	|\sin z|^2 \ge \sinh^2 y, \quad |\cos z|^2 \ge \sinh^2 y,
\end{equation}
\begin{equation}\label{2.7}
	y{\partial \over \partial y }|\sin z|^2= 
	y{\partial \over \partial y}|\cos z|^2 \ge 2 y^2,
\end{equation}
\begin{equation}\label{2.8}
	{\partial^2 \over \partial y^2}|\sin z|^2 =
	{\partial^2 \over \partial y^2}|\cos z|^2 \ge 2 \cosh^2 y 
\end{equation}
\begin{equation}\label{2.9}
	{\partial^2 \over \partial y^2}|\sin z|^2 =
	{\partial^2 \over \partial y^2}|\cos z|^2 \ge 2 + 2 \sinh^2 y 
\end{equation}
\noindent
which are consequences of (2.1)--(2.4).

\section{Bessel functions and ${}_0F_1(c\/; z)$ functions}
\setcounter{equation}{0}

Since the identities  and inequalities in \S 2 give the reality of the 
zeros of the Bessel functions \cite[(1.71.2)]{sz}
\begin{equation}\label{3.1}
	J_{-{1\over2}}(z)=({2 \over{\pi z}})^{1 \over 2}\cos z, 
	\quad J_{1\over2}(z)=({2 \over {\pi z}})^{1 \over 2}\sin z,
\end{equation}
this suggests that it should be possible to use sums of squares
to prove Lommel's theorem (see Watson \cite[p.~482]{wat}) that all of 
the zeros of the Bessel function \cite[7.2(3)]{erd}
\begin{equation}\label{3.2}
	J_\alpha(z)={(z/2)^{\alpha} \over \Gamma(\alpha+1)}\;
	{}_0F_1(\alpha+1;-z^2/4)
\end{equation}
\noindent
are real when $\alpha>-1$.  With this aim in mind and in order to work 
with entire functions, we set
\begin{equation}\label{3.3}
	{\cal J}_\alpha(z)=z^{-\alpha}J_\alpha(z)={2^{-\alpha} \over
	\Gamma(\alpha+1)}\; {}_0F_1(\alpha+1;-z^2/4),
\end{equation}
which is an even entire function of $z$ such that
${\overline{{\cal J}_\alpha(z)} = {\cal J}_\alpha(\overline z)}$
when $\alpha$ is real.  

Let $\alpha>-1$. Then, from the product formula (37) in Carlitz \cite{car},
\begin{equation}\label{3.4}
	|{\cal J}_\alpha(z)|^2 = \sum_{k=0}^\infty
	\frac{(\alpha+\frac12)_k2^{k-\alpha}}{k!(2\alpha+1)_k\Gamma(\alpha+1)}
	(x^2+y^2)^k{\cal J}_{\alpha+k}(2x).
\end{equation}
To express each of the Bessel functions on the right side of (3.4) as
a sum of squares of real valued Bessel functions 
observe that from the addition theorem for
Bessel functions \cite[7.15(30)]{erd} we have the expansion
\begin{equation}\label{3.5}
	{\cal J}_{\alpha+k}(2x) = 2^{k+\alpha}\Gamma(k+\alpha)\sum_{j=0}^\infty
	\frac{(j+k+\alpha)(2k+2\alpha)_j}{j!} (-1)^j x^{2j}
	({\cal J}_{\alpha+j+k}(x))^2.
\end{equation}
Hence, substituting (3.5) into (3.4) and changing the order of summation
we find that
\begin{multline}\label{3.6}
	|{\cal J}_\alpha(z)|^2 
	=\sum_{n=0}^\infty
	{(n+\alpha)(2\alpha)_n \over{\alpha \; n!}} (-1)^n x^{2n} \\
	\times {}_2F_1(-n,n+2\alpha\/;2\alpha+1;1+y^2/x^2) 
	({\cal J}_{\alpha+n}(x))^2.
\end{multline}
Now apply the Euler transformation formula \cite[2.9(3)]{erd}
\begin{equation}\label{3.7}
	{}_2 F_1(a,b\/;c\/;z) = (1-z)^{-a} {}_2 F_1(a,c-b\/;c\/;z/(z-1))
\end{equation}
to the above ${}_2F_1$ series to obtain the desired sum of squares
expansion formula
\begin{align}\label{3.8}
	|{\cal J}_\alpha(z)|^2 &= ({\cal J}_{\alpha}(x))^2
	+2(\alpha+1)  y^2 ({\cal J}_{\alpha+1}(x))^2 \\ \notag
	&+ \sum_{n=2}^\infty{(2n+2\alpha)(2\alpha+1)_{n-1} \over{ n!}} y^{2n}\\
	&\times {}_2F_1(-n,1-n\/;2\alpha+1;1+x^2/y^2) 
	({\cal J}_{\alpha+n}(x))^2. \notag
\end{align}
When $n \ge 2, \alpha>-1$ and $y \ne 0,$ the positivity of the
coefficients of $({\cal J}_{\alpha+n}(x))^2$ in (3.8) follows from
\begin{multline}\label{3.9}
	(2\alpha+1)\> {}_2F_1(-n,1-n\/; 2\alpha+1; 1+x^2/y^2)
	=(2\alpha+1+n^2-n) \\
	+ n(n-1)x^2/y^2 + \sum^n_{k=2}
	\frac{(-n)_k(1-n)_k}{k!\,(2\alpha+2)_{k-1}} (1+x^2/y^2)^k > 0.
\end{multline}
Hence, since the real zeros of ${\cal J}_\alpha(x)$ and ${\cal
J}_{\alpha+1}(x)$ are interlaced, (3.8) gives a sum of squares proof
that the entire functions ${\cal J}_\alpha(z)$, and thus the Bessel
functions $ J_\alpha(z),$ have only real zeros when $\alpha>-1$. Letting
$\alpha \to -1$ it follows that the Bessel function
${\cal J}_{-1}(z) = \lim_{\alpha\to-1}$
${\cal J}_{\alpha}(z) = -{\cal J}_1 (z)$ has only real zeros.

Notice that the inequality
\begin{equation}\label{3.10}
	|{\cal J}_\alpha(z)|^2 \ge ({\cal J}_{\alpha}(x))^2
	+2(\alpha+1)  (y{\cal J}_{\alpha+1}(x))^2 >0, \hfill
	\>	y \ne0, \> \alpha > -1,
\end{equation}
and in fact infinitely many inequalities follow from (3.8) by just
dropping terms from the right side of (3.8). Analogous to (2.7)--(2.9),
it follows by differentiating equation (3.6) with respect to $y$
and applying (3.7) that we also have the identities
\begin{multline}\label{3.11}
	y {\partial\over \partial y}|{\cal J}_\alpha(z)|^2 
	= 4y^2 \sum_{n=0}^\infty
	{(n+\alpha+1)(2\alpha+2)_n \over{ n!}} y^{2n} \\
	\times {}_2F_1(-n,-n\/;2\alpha+2;1+x^2/y^2) 
	({\cal J}_{\alpha+n+1}(x))^2 
\end{multline}
and
\begin{eqnarray}\label{3.12} \hskip3em
	\frac{\partial^2}{\partial y^2}|{\cal J}_\alpha(z)|^2
	&=& 4 \sum_{n=0}^\infty\frac{(n+\alpha+1)(2\alpha+2)_n}{n!}
	y^{2n} \\
	&\times& {}_2F_1(-n,-n\/;2\alpha+2;1+x^2/y^2) 
	({\cal J}_{\alpha+n+1}(x))^2 \notag \\ [2mm]
	&+&8y^2 \sum_{n=0}^\infty
	{(n+\alpha+2)(2\alpha+3)_{n+1} \over{ n!}} y^{2n} \notag \\
	&\times& {}_2F_1(-n,-n-1;2\alpha+3;1+x^2/y^2) 
	({\cal J}_{\alpha+n+2}(x))^2, \notag
\end{eqnarray}
which give infinitely many inequalities, such as, e.g.,
\begin{equation}\label{3.13}
	y{\partial\over \partial y}|{\cal J}_\alpha(z)|^2\ge 4(\alpha+1)
	(y{\cal J}_{\alpha+1}(x))^2 \ge 0, \quad \alpha\ge-1,
\end{equation}
\begin{equation}\label{3.14}
	{\partial^2\over \partial y^2}|{\cal J}_\alpha(z)|^2\ge
	4(\alpha+1)({\cal J}_{\alpha+1}(x))^2 \ge 0, \quad \alpha\ge-1,
\end{equation}
each of which proves that $J_\alpha(z)$ has only real zeros when
$\alpha\ge-1.$

In view of (3.3) the reality of the zeros of ${\cal J}_\alpha(z)$ when
$\alpha > -1$ is equivalent to the statement that all of the zeros of
the confluent hypergeometric function ${}_0F_1(c\/;z)$ are real and
negative when $c>0.$  However, it is known \cite{hil} that the zeros of
${}_0F_1(c\/; z)$      
are also real (but not necessarily negative) when
$-1<c<0.$  Because this fact does not follow from (3.8) or (3.11)--(3.14),
 we will next show how it can also be proved by using sums of
squares of real value functions.

 From formulas (53) and (52) in Burchnall and Chaundy \cite{bc41} it follows
that if $c$ is real valued and $c \ne 0, -1, -2, \dots,$ then we have
the expansion formulas
\begin{equation}\label{3.15}
	\left|{}_0F_1(c\/;z)\right|^2 = \sum^\infty_{k=0}
	\frac{1}{k!\,(c)_k(c)_{2k}} (x^2+y^2)^k {}_0F_1(c+2k\/;2x)
\end{equation}
and
\begin{equation}\label{3.16}
	{}_0F_1(c+2k\/;2x) = \sum^\infty_{j=0} \frac{(-1)^j}
	{j!\, (c+2k+j-1)_j(c+2k)_{2j}} x^{2j}\left({}_0F_1(c+2k+2j\/;x)\right)^2.
\end{equation}
As in the Bessel function case, substitute (3.16) into (3.15) and
change the order of summation to get
\begin{multline}\label{3.17}
	\left|{}_0F_1(c\/;z)\right|^2 = \sum^\infty_{n=0}  \frac{(-1)^n}
	{n!\, (c+n-1)_n (c)_{2n}} x^{2n} \\
	\times {}_2F_1(-n,n+c-1; c\/; 1+y^2/x^2)({}_0F_1(c+2n\/;x))^2
\end{multline}
which, by applying the transformation formula (3.7), gives
\begin{multline}\label{3.18}
	\left|{}_0F_1(c\/;z)\right|^2 = \sum^\infty_{n=0}  \frac{1}
	{n!\,(n+c-1)_n(c)_{2n}} y^{2n} \\
	\times {}_2F_1(-n,1-n\,;c\/;1+x^2/y^2)({}_0F_1(c+2n\/;x))^2.
\end{multline}
When $c>0$ and $y \ne 0$  the coefficient of $({}_0F_1(c+2n\/; x))^2$ in 
the series in (3.18) is obviously positive.
Hence, since ${}_0F_1(c\/; x) > 0$ when $c>0$ and $x\ge 0,$ 
(3.18) gives another proof that  ${}_0F_1(c\/; z)$ has
only real negative zeros when $c>0.$  

To handle the case $-1<c<0$  differentiate equation (3.17) with respect to 
$y$ and apply (3.7) to obtain
\begin{multline}\label{3.19}
	y\frac{\partial}{\partial y}\left|c\>{}_0F_1(c\/;z)\right|^2 = 2y^2
	\sum^\infty_{n=0} \frac{(c+1)_n}{n!\,(c+1)_{2n}(c+1)_{2n+1}}y^{2n} \\
	\times {}_2F_1(-n,-n; c+1; 1+x^2/y^2)({}_0F_1(c+2n+2;x))^2
\end{multline}
and
\begin{align}\label{3.20}
	\frac{\partial^2}{ \partial y^2}\left|c\>{}_0F_1(c\/;z)\right|^2
	  &= 2\sum^\infty_{n=0}
	  \frac{(c+1)_n}{n!\,(c+1)_{2n}(c+1)_{2n+1}}y^{2n} \\ \notag
	&\times {}_2F_1(-n,-n\/;c+1;1+x^2/y^2)({}_0F_1(c+2n+2;x))^2 \\ \notag
	&+ 4y^2 \sum^\infty_{n=0}
	  \frac{(c+2)_{n+1}}{n!\,(c+1)_{2n+2}(c+1)_{2n+3}} y^{2n} \\ \notag
	&\times {}_2F_1(-n,-n-1;c+2;1+x^2/y^2)({}_0F_1(c+2n+4;x))^2
\end{align}
which, in particular, give the inequalities
\begin{equation}\label{3.21}
	y\frac{\partial}{\partial y}\big|c(c+1)\>{}_0F_1(c\/;z)\big|^2 \ge
	2(c+1)(y\>{}_0F_1(c+2;x))^2,\quad c\ge-1,
\end{equation}
and
\begin{equation}\label{3.22} \quad
    \frac{\partial^2}{\partial y^2}\left|c(c+1)\>{}_0F_1(c;z)\right|^2 \ge 
	2(c+1)({}_0F_1(c+2;x))^2,\quad c\ge-1.
\end{equation}
Since the coefficients on the right hand sides of (3.19)--(3.22) are
clearly positive when $c > -1$ and $y \ne 0,$ these formulas prove that
the functions $c(c+1) \> {}_0F_1(c\/; z)$ have only real zeros when $c\ge
-1,$ where it is understood that $c(c+1) \> {}_0F_1(c\/; z)$ is to be
replaced by its $c \rightarrow 0$ limit case $z \> {}_0F_1(2; z)$ when
$c=0,$ and by its $c \rightarrow -1$  limit case $z^2 \> {}_0F_1(3;
z)/2$ when $c=-1.$

\section{Laguerre polynomials and ${}_1F_1(a\/; c\/; z)$ functions}
\setcounter{equation}{0}

When $\alpha > - 1$ the Laguerre polynomials
\begin{equation}\label{4.1}
	L_n^\alpha(z) = \frac{(\alpha+1)_n}{n!}\>{}_1F_1(-n\/;\alpha+1; z)
\end{equation}
satisfy the orthogonality relation
\begin{equation}\label{4.2}
	\int^\infty_0 L_n^\alpha(x) L_m^\alpha(x) x^\alpha e^{-x}\, dx =
	\frac{\Gamma(n+\alpha+1)}{n!} \delta_{nm},\;\; n,m=0,1,2,\dots,
\end{equation}
from which it follows by a standard argument (cf. \cite[\S3.3]{sz})
that the zeros of  $L_n^\alpha(z)$ are real and positive. Analogous to
the last part of the previous section, in this section we will derive
some sums of squares expansions which, besides proving the reality of
the zeros of these polynomials when $\alpha > - 1,$ also prove that
they have only real zeros (not necessarily positive) when $-1\ge\alpha
\ge -2,$ where $L_n^\alpha(z)$ is defined to be the $\alpha \rightarrow
-k$ limit case of (4.1) when $\alpha$ is a negative integer $-k.$ Thus
$L_1^\alpha(z) = \alpha + 1 - z,$ which has a negative zero when
$\alpha < -1,$ and $L_2^\alpha(z) = ((\alpha+1)(\alpha+2)-
2(\alpha+2)z+z^2)\big/2,$ which has non-real zeros when $\alpha < -2.$

Let $\alpha$ be real valued. Substituting the sum of squares of Laguerre
polynomials expansion (from \cite[(91)]{bc41})
\begin{align}\label{4.3}
	L^{\alpha+2k}_{n-k} (2x) &= \sum^{n-k}_{j=0}
	\frac{(n-k-j)!\, (2k+2j+\alpha)(2k+\alpha)_j}
	  {j!\, (2k+\alpha)(2k+\alpha+1)_{n+j-k}} \\ \notag	
	&\times (-1)^j x^{2j} \left( L^{\alpha+2k+2j}_{n-k-j} (x)\right)^2
\end{align}
into the special case of \cite[(5.4)]{ba}
\begin{equation}\label{4.4}
	\left| L_n^\alpha(z)\right|^2 = \frac{(\alpha+1)_n}{n!} 
	\sum^n_{k=0}\frac{1}{k!\,(\alpha+1)_k} (x^2+y^2)^k
	L^{\alpha+2k}_{n-k} (2x)
\end{equation}
and changing the order of summation yields
\begin{multline}\label{4.5}
    \left| L_n^\alpha(z)\right|^2 = \frac{(\alpha+1)_n}{n!}\sum^n_{k=0}
	\frac{(n-k)!\,(2k+\alpha)(\alpha)_k}{k!\,\alpha(\alpha+1)_{n+k}}
	(-1)^k x^{2k} \\
	\times {}_2F_1(-k,k+\alpha\/;\alpha+1;1+y^2/x^2)
	\left(L^{\alpha+2k}_{n-k} (x)\right)^2.
\end{multline}
Then application of (3.7) gives
\begin{multline}\label{4.6}
    \left| L_n^\alpha(z)\right|^2 = \frac{(\alpha+1)_n}{n!}\sum^n_{k=0}
    \frac{(n-k)!\,(2k+\alpha)(\alpha)_k}{k!\,\alpha(\alpha+1)_{n+k}}
	y^{2k} \\
	\times {}_2F_1(-k,1-k\/;\alpha+1;1+x^2/y^2)
	\left(L^{\alpha+2k}_{n-k} (x)\right)^2.
\end{multline}
Since $L_0^\alpha (x) \equiv 1$ and the coefficients on the right hand side
 of (4.6) are clearly positive when $\alpha > -1$ and $y \ne 0,$
the expansion (4.6) proves that the Laguerre polynomials have only real zeros
 when $\alpha > -1.$  This also follows, in particular, from the inequalities
\begin{equation}\label{4.7}
    \left| L_n^\alpha(z)\right|^2 \ge
	\frac{(\alpha+1)_n}{n!\, n!\,(n+\alpha)_n} y^{2n}
	\>{}_2F_1(-n, 1-n\/;\alpha+1;1+x^2/y^2),  \quad  \alpha>-1,
\end{equation}
and 
\begin{equation}\label{4.8}
	\left| L_n^\alpha(z)\right|^2 \ge \left| L_n^\alpha(x)\right|^2 +
	\frac{(\alpha+1)_n}{n!\, n!\,(n+\alpha)_n} y^{2n},
	\quad  \alpha>-1,\ n\ge1,
\end{equation}
which are consequences of (4.6).

Now differentiate equation (4.5) with respect to $y$ and apply (3.7) to
derive the expansions
\begin{multline}\label{4.9}
	y\frac{\partial}{\partial y} \left| L_n^\alpha(z)\right|^2 =
	2y^2 \sum^{n-1}_{k=0}
	\frac{(n-k-1)!\,(2k+\alpha+2)(\alpha+2)_k}{n!\,k!\,(n+\alpha+1)_{k+1}}
	y^{2k} \\
	\times {}_2F_1(-k,-k\/;\alpha+2;1+x^2/y^2)
	\left( L^{\alpha+2k+2}_{n-k-1} (x)\right)^2, \quad n\ge 1,
\end{multline}
and
\begin{align}\label{4.10}
	\frac{\partial^2}{\partial y^2} \left| L_n^\alpha(z)\right|^2 &=
	  2 \sum^{n-1}_{k=0}
	  \frac{(n-k-1)!\,(2k+\alpha+2)(\alpha+2)_k}{n!\,k!\,(n+\alpha+1)_{k+1}}
	  \>y^{2k} \\ \notag
	&\times {}_2F_1(-k,-k\/;\alpha+2;1+x^2/y^2)
	  \left( L^{\alpha+2k+2}_{n-k-1} (x)\right)^2 \\[2mm] \notag
	&+ 4y^2 \sum^{n-2}_{k=0} 
	  \frac{(n-k-2)!\,(2k+\alpha+4)(\alpha+3)_{k+1}}{n!\,k!\,(n+\alpha+1)_{k+2}}
	  \>y^{2k} \\ \notag
	&\times {}_2F_1(-k,-k-1;\alpha+3;1+x^2/y^2)
	  \left( L^{\alpha+2k+4}_{n-k-2} (x)\right)^2\!, \quad n\ge1,
\end{align}
which yield, e.g., the inequalities
\begin{equation}\label{4.11}\quad
	 y\frac{\partial}{\partial y} \left| L_n^\alpha(z)\right|^2 \ge
	 \frac{2(\alpha+2)}{n(n+\alpha+1)}
	 \left(yL^{\alpha+2}_{n-1}(x)\right)^2, \quad \alpha>-2,\ n\ge1,
\end{equation}
and
\begin{equation}\label{4.12}\quad
	 \frac{\partial^2}{\partial y^2} \left| L_n^\alpha(z)\right|^2 \ge
	 \frac{2(\alpha+2)}{n(n+\alpha+1)}
	 \left(L^{\alpha+2}_{n-1}(x)\right)^2, \quad \alpha>-2,\ n\ge1,
\end{equation}
and prove (after letting $\alpha\to-2$) that the polynomials
$L_n^\alpha (z)$ have only real zeros when $\alpha \ge -2.$

For the confluent hypergeometric functions 
${}_1F_1(a\/; c\/;z)$ with $a$ and $c$
real valued and $c \ne 0, -1, -2, \dots,$
use of the expansion formulas \cite[(42) and (43)]{bc41}
instead of (4.3) and (4.4)
yields the nonterminating extension of (4.5)
\begin{multline}\label{4.13}
	|{}_1F_1(a\/; c\/;z)|^2 = \sum^\infty_{k=0}\frac{(a)_k(c-a)_k} 
	{k!\,(c)_{2k}(c+k-1)_k} x^{2k} \\
	\times {}_2F_1(-k, c+k-1; c\/;1+y^2/x^2)({}_1F_1(a+k\/;c+2k\/;x))^2
\end{multline}
and hence, by (3.7),
\begin{multline}\label{4.14}
	|{}_1F_1(a\/; c\/;z)|^2 = \sum^\infty_{k=0}\frac{(a)_k(c-a)_k}
	{k!\,(c)_{2k}(c+k-1)_k} (-1)^k y^{2k}  \\
	\times {}_2F_1(-k,1-k\/;c\/;1+x^2/y^2)({}_1F_1(a+k\/;c+2k\/;x))^2.
\end{multline}
Then differentiation of equation (4.13) with respect of $y$ and application of 
(3.7) gives the following extensions of (4.9) and (4.10)   
(and also of (3.19) and (3.20)), respectively,
\begin{multline}\label{4.15}
	y\frac{\partial}{\partial y}|c(c+1)\>{}_1F_1(a\/; c\/;z)|^2 =
	2y^2\sum^\infty_{k=0}\frac{(a)_{k+1}(c-a)_{k+1}(c+1)}
	{k!\,(c+2)_{2k}(c+k+1)_k} (-1)^{k+1} y^{2k}  \\
	\times {}_2F_1(-k,-k\/;c+1;1+x^2/y^2)({}_1F_1(a+k+1;c+2k+2;x))^2
\end{multline}
and
\begin{align}\label{4.16}
	\frac{\partial^2}{\partial y^2} \, &|c(c+1)\>{}_1F_1(a\/;c\/;z)|^2 
		= 2\sum^\infty_{k=0} \frac{(a)_{k+1}(c-a)_{k+1}(c+1)}
		{k!\,(c+2)_{2k}(c+k+1)_k} (-1)^{k+1} y^{2k} \\ \notag
	&\times {}_2F_1(-k,-k; c+1;1+x^2/y^2)({}_1F_1(a+k+1;c+2k+2;x))^2 \\ \notag
	&+ 4y^2\sum^\infty_{k=0}\frac{(a)_{k+2}(c-a)_{k+2}}
		{k!\,(c+2)_{2k+2}(c+k+3)_k}(-1)^ky^{2k} \\ \notag
	&\times {}_2F_1(-k,-k-1;c+2;1+x^2/y^2)({}_1F_1(a+k+2;c+2k+4;x))^2.
\end{align}

If $a = -n$ is a negative integer and $c = \alpha +1,$ then (4.13)--(4.16)
 reduce to (4.5), (4.6), (4.9), (4.10), respectively. If $a = c+n$ with $n$
 a nonnegative integer, then (4.15) and (4.16) reduce to terminating
sums of squares expansions with nonnegative coefficients which prove that 
$c(c+1) \>{}_1F_1(c+n\/; c\/;z),$ as a 
function of $z,$ has only real zeros when $c \ge -1,$ where this function
is to be replaced by its  $c \rightarrow 0$ and $c \rightarrow -1$  
limit cases when  $c=0$ and $c=-1,$ respectively.  
It should be noted that, in view of Kummer's transformation formula 
\cite[6.3(7)]{erd}
\begin{equation}\label{4.17}
	{}_1F_1(a\/;c\/;x) = e^x\>{}_1F_1(c-a\/;c\/;-x),
\end{equation}
these results on the zeros of $c(c+1)\> {}_1F_1(c+n\/; c\/;z)$
are equivalent to those obtained above for the Laguerre polynomials.

\section{Jacobi polynomials}
\setcounter{equation}{0}

When $ \alpha > -1$ and $\beta > -1$ the Jacobi polynomials
\begin{equation}\label{5.1}
	P^{(\alpha,\beta)}_n (z) = \frac{(\alpha+1)_n}{n!}\>
	{}_2F_1(-n,n+\alpha+\beta+1;\alpha+1;(1-z)/2)
\end{equation}
satisfy the orthogonality relation
\begin{equation}\label{5.2}
	\int^1_{-1} P^{(\alpha,\beta)}_n (x)P^{(\alpha,\beta)}_m(x)(1-x)^\alpha
	(1+x)^\beta\>dx = 0,\quad n\ne m,
\end{equation}
for $n, m = 0, 1, 2, \dots,$ and hence, by \cite[Theorem 3.3.1]{sz}, these
polynomials have only real zeros. In our derivation of sums of squares 
expansions which imply the reality of the zeros of these
polynomials we will start out by deriving sums of squares 
expansions for nonterminating ${}_2F_1(a, b\/; c\/; z)$ hypergeometric series 
with $|z| < 1$ (for convergence).

Let $a, b, c$ be real valued, $c \ne 0, -1, -2, \dots,$ and $|z| < 1.$ 
Then formula \cite[(51)]{bc40} gives the expansion
\begin{align}\label{5.3}
	|{}_2F_1(a,b\/;c\/;z)|^2 &= \sum^\infty_{k=0}
	\frac{(a)_k(b)_k(c-a)_k(c-b)_k}{k!\,(c)_k(c)_{2k}} (x^2+y^2)^k \\
	&\times {}_2F_1(a+k,b+k\/;c+2k\/;2x-x^2-y^2). \notag
\end{align}
Unfortunately, application of the inversion \cite[(50)]{bc40} of
\cite[(51)]{bc40} to each of the ${}_2F_1(a+k, b+k\/; c+2k\/;
2x-x^2-y^2)$ functions on the right side of equation (5.3) just returns
one back to the function that is on the left side. Therefore, we use
formulas (44), (45), (50) in \cite{bc40} to obtain, respectively, the
expansions
\begin{multline}\label{5.4}
	{}_2F_1(a+k,b+k\/;c+2k\/;2x-x^2-y^2)  \\
	  = \sum^\infty_{j=0}\frac{(a+k)_j(b+k)_j}{j!\,(c+2k)_j}(-1)^j
	  (x^2+y^2)^j\> {}_2F_1(a+k+j,b+k+j\/;c+2k+j\/;2x),
\end{multline}
\begin{multline}\label{5.5}
	{}_2F_1(a+k+j,b+k+j\/;c+2k+j\/;2x) 
	  = \sum^\infty_{m=0}\frac{(a+k+j)_m(b+k+j)_m}{m!\,(c+2k+j)_m} x^{2m} \\
	  \times {}_2F_1(a+k+j+m,b+k+j+m\/;c+2k+j+m\/;2x-x^2),
\end{multline}
\begin{multline}\label{5.6}
	{}_2F_1(a+k+j+m,b+k+j+m\/;c+2k+j+m\/;2x-x^2)  \\
	  = \sum^\infty_{n=0}\frac{(a+k+j+m)_n(b+k+j+m)_n(c-a+k)_n(c-b+k)_n}
	  {n!\,(c+2k+j+m+n-1)_n(c+2k+j+m)_{2n}}  (-1)^n x^{2n}  \\
	  \times({}_2F_1(a+k+j+m+n,b+k+j+m+n\/;c+2k+j+m+2n\/;x))^2,
\end{multline}
and then substitute these expansions in turn into (5.3), change the order
of summation and use the binomial theorem to obtain 
\begin{multline}\label{5.7}
	|{}_2F_1(a,b\/;c\/;z)|^2
	  = \sum^\infty_{m=0}\sum^m_{j=0}\frac{(a)_m(b)_m(c-a)_j(c-b)_j}
	  {j!\, (m-j)!\,(c)_{m+j}(m+c-1)_j} (-1)^m x^{2j}y^{2m-2j}  \\
	  \times{}_2F_1(-j,m+c-1;c\/;1+y^2/x^2)({}_2F_1(m+a,m+b\/;m+j+c\/;x))^2.
\end{multline}

Application of (3.7) to the first ${}_2F_1$ series on the right
side of (5.7) gives
\begin{multline}\label{5.8}
	|{}_2F_1(a,b\/;c\/;z)|^2 
	  = \sum^\infty_{m=0}\sum^m_{j=0}\frac{(a)_m(b)_m(c-a)_j(c-b)_j}
	  {j!\, (m-j)!\,(c)_{m+j}(m+c-1)_j} (-1)^{m+j} y^{2m}  \\
	  \times{}_2F_1(-j,1-m\/;c\/;1+x^2/y^2)({}_2F_1(m+a,m+b\/;m+j+c\/;x))^2,
\end{multline}
which contains (4.14) as a limit case.  When $a = -n$ is a negative
integer, $b=n+\alpha+\beta+1$ and $c=\alpha+1,$ it follows from (5.8) that
\begin{multline}\label{5.9}
	\left|\frac{n!}{(\alpha+1)_n} P^{(\alpha,\beta)}_n (1-2z)\right|^2 \\
	  = \sum^n_{m=0}\sum^m_{j=0}\frac{(-n)_m(n+\alpha+\beta+1)_m
	  (n+\alpha+1)_j(-n-\beta)_j}{j!\,(m-j)!\,(\alpha+1)_{m+j}(m+\alpha)_j}
	  (-1)^{m+j} y^{2m} \\
	  \times{}_2F_1(-j,1-m\/;\alpha+1;1+x^2/y^2)
	  ({}_2F_1(m-n,m+n+\alpha+\beta+1;m+j+\alpha+1;x))^2,
\end{multline}
which gives a sums of squares proof that the Jacobi polynomials 
$P_n^{(\alpha, \beta)}(z)$ have only real zeros when $ \alpha, \beta > -1$
(since the coefficients in (5.9) are then clearly positive) and hence, by
continuity, when  $ \alpha, \beta \ge -1.$  The restriction that $ \alpha, 
\beta \ge -1$ cannot be extended to $ \alpha, \beta \ge -2$
because $P_2^{(\alpha, \beta)}(z)$
has non-real zeros when $ \alpha, \beta > -2$ and $\alpha+\beta<-3.$

As in sections 3 and 4 one may repeatedly differentiate (5.7)
with respect to $y$ and apply (3.7) to obtain extensions of (4.15), (4.16), etc.
But, since the resulting identities are quite lengthly and do not add 
any additional $(\alpha, \beta)$ for which the Jacobi polynomials have only 
real zeros, we will omit them and only point out that the first two 
differentiations give identities that yield, in particular, the inequalities
\begin{equation}\label{5.10}
	y\frac{\partial}{\partial y}\left|P^{(\alpha,\beta)}_n (1-2z)\right|^2
	  \ge \frac{2n(n+\alpha+\beta+1)_n(\alpha+1)_n}{n!\,n!} y^{2n}
\end{equation}
and
\begin{equation}\label{5.11}
	\frac{\partial^2}{\partial y^2}\left|P^{(\alpha,\beta)}_n (1-2z)\right|^2
	  \ge \frac{2n(2n-1)(n+\alpha+\beta+1)_n(\alpha+1)_n}{n!\,n!} y^{2n-2}
\end{equation}
when $n\ge 1$ and $ \alpha, \beta \ge -1.$

In subsequent papers it will be shown that squares of real valued functions
can also be used to prove the reality of the zeros of some non-classical
families of orthogonal polynomials, of the cosine transforms
$$ \int_0^\infty  e^{-a \cosh t} \cos zt \, dt, \qquad a > 0,$$
and of some other entire functions.


George Gasper

Department of Mathematics

Northwestern University

Evanston, IL\ \ 60208

E-Mail: g-gasper\char`\@ nwu.edu

\end{document}